\documentclass{amsart}

\usepackage{amssymb}
\usepackage[all]{xy}
\usepackage{hyperref}
\usepackage{mathrsfs}
\usepackage{amsthm}
\newtheorem{thm}{Theorem}%[section]

\newtheorem{lem}[thm]{Lemma}
\newtheorem{cor}[thm]{Corollary}

\newtheorem{prop}[thm]{Proposition}

\theoremstyle{definition}

\newtheorem{say}[thm]{}
\newtheorem{exmp}[thm]{Example}

\newtheorem{rem}[thm]{Remark}          

\newtheorem*{ack}{Acknowledgments}      % \renewcommand{\theack}{} 
\newtheorem{notation}[thm]{Notation}   
  
\newtheorem{defn-thm}[thm]{Definition--Theorem}  %!!!!!!!!!!!!!!!!!!!!!!!!
\newtheorem{defn-lem}[thm]{Definition--Lemma}  %!!!!!!!!!!!!!!!!!!!!!!!!
\theoremstyle{remark}

\let \cedilla =\c
\renewcommand{\c}[0]{{\mathbb C}}  

\renewcommand{\o}[0]{{\mathcal O}} 
\newcommand{\z}[0]{{\mathbb Z}}

\renewcommand{\r}[0]{{\mathbb R}} 

\renewcommand{\a}[0]{{\mathbb A}}

\newcommand{\p}[0]{{\mathbb P}}

\newcommand{\qtq}[1]{\quad\mbox{#1}\quad}
\newcommand{\spec}[0]{\operatorname{Spec}}

\newcommand{\aut}[0]{\operatorname{Aut}}

\newcommand{\chr}[0]{\operatorname{char}}

\newcommand{\grass}[0]{\operatorname{Gr}}

\newcommand{\onto}[0]{\twoheadrightarrow}

\newcommand{\an}[0]{\operatorname{an}}

\def\loccoh#1.#2.#3.#4.{H^{#1}_{#2}(#3,#4)}

\DeclareMathAlphabet{\mathchanc}{OT1}{pzc}%
                                {m}{it}

\newcommand{\gm}[0]{{\mathbb G}_m}
\newcommand{\ga}[0]{{\mathbb G}_a}

\newcommand{\GL}{\mathrm{GL}}
\newcommand{\PGL}{\mathrm{PGL}}

\newcommand{\SL}{\mathrm{SL}}

\newcommand{\GO}{\mathrm{GO}}
\newcommand{\OO}{\mathrm{O}}

\newcommand{\Sp}{\mathrm{Sp}}

\newcommand{\sym}[0]{\operatorname{Sym}}

\usepackage[all]{xy}\xyoption{dvips}

\begin{document}
\bibliographystyle{amsalpha}

%\hfill\today

 \title[Fulton-Hansen theorem]{A Fulton-Hansen theorem for \\ almost homogeneous spaces}
\author{J\'anos Koll\'ar and Aaron Landesman}
\dedicatory{Dedicated to Fabrizio Catanese on his 70th birthday}
\begin{abstract} We prove a generalization of the Fulton-Hansen connectedness theorem, where ${\mathbb P}^n$ is replaced by a normal variety on which an algebraic group acts with a dense orbit.
\end{abstract}

 \maketitle

\bigskip

The  product version of the Fulton-Hansen theorem  \cite{fultonH:a-connectedness-theorem}
says that
if $Y_1, Y_2$ are irreducible, projective varieties, 
  $Y_i\to \p^n$ are finite morphisms
and $\dim Y_1+\dim Y_2>\dim \p^n$, then   $Y_1\times_{\p^n}Y_2$ is connected.
As Catanese emphasized both in the panoramic survey
\cite{catanese:topological-methods-in-moduli-theory}
and
in the shorter version \cite{catanese:topological-methods-in-algebraic-geometry},
one should  look at the topological side of every question in algebraic geometry.  Thus one should view the  Fulton-Hansen theorem as a
generalization of Lefschetz's hyperplane section theorem for the fundamental group.
In this form, the  quasi-projective version 
says that
if $Y_1, Y_2\subset \p^n$ are normal, irreducible, quasi-projective varieties, 
$p_i:Y_i\to \p^n$ are quasi-finite morphisms, and $\dim Y_1+\dim Y_2>\dim \p^n$, then 
$$
\pi_1\bigl(Y_1\times_{\p^n, g_1, g_2}Y_2\bigr)\onto \pi_1(Y_1)\times \pi_1(Y_2) \qtq{is surjective}
$$
for general $g_1, g_2\in \aut(\p^n)$; here our notation means that we first
compose $p_i$ with $g_i$ and then take fiber product. 
See  
\cite{fultonH:a-connectedness-theorem, deligne:le-groupe-fondamental}
and
\cite[Remark 9.3]{fultonL:connectivity-and-its-applications} for the original
proofs and Corollaries~\ref{ful-han.thm.cor}--\ref{ful-han.thm.cor-ii}  for further discussion.
Also see \cite[Theorem 3.3.6]{lazarsfeld:positivity-i} and
\cite[Theoreme 7.2]{jouanolou1982theoremes}.

%% As a minor technical point, the original Fulton-Hansen theorem begins with a map
%% $Z \to \p^n \times \p^n$ and shows the preimage of the diagonal is connected. In
%% the above discussion, we assumed $Z = Y_1 \times Y_2$. However, we also recover
%% the Fulton-Hansen theorem in this greater generality, at least when the above map is
%% finite, in Corollary~\ref{ful-han.thm.cor-ii}.
%% Still, the special case that $Z = Y_1 \times Y_2$ gives many interesting
%% applications, and is often useful to keep in mind.

Our main result is Theorem~\ref{dense.orb.thm}, which generalizes the above
results.
We replace 
$(\PGL_{n+1}, \p^n)$ by a pair $(G, X)$, where $X$ is a normal quasi-projective variety
and  $G$ is an  algebraic group acting on  $X$ with a fixed point and a dense open orbit. 
We obtain the Fulton-Hansen theorem by applying
Theorem~\ref{dense.orb.thm} to  $X=\p^n$ and $G\subset \PGL_{n+1}$ the
stabilizer of a  point; see Corollary~\ref{ful-han.thm.cor}.
This also implies the original form of the Fulton-Hansen theorem, which begins with a map
$Z \to \p^n \times \p^n$ and shows that the preimage of the diagonal is connected, see Corollary~\ref{ful-han.thm.cor-ii}.

The usual approaches to Fulton-Hansen--type  theorems exploit ampleness
properties  of the tangent bundle  $T_X$. We get the strongest result when $T_X$
is ample, which holds only for $X=\p^n$
\cite{mori:projective-manifolds-with-ample-tangent-bundles}.
For other homogeneous spaces  $X=G/P$, one can develop a  theory that measures the failure of ampleness    \cite{sommese:submanifolds-of-abelian-varieties, sommese:complex-subspaces-i, fultonL:connectivity-and-its-applications, faltings:formale-geometrie-und-homogene-raume, goldstein:ampleness-and-connectedness}.  
This leads to results that give  $\pi_1$-surjectivity if
$\dim Y_1+\dim Y_2$ is large enough. 
For example, let $0\leq r\leq n/2$ and let $X=\grass(\p^r, \p^n)$.
By \cite[Connectedness Theorem, p. 361]{goldstein:ampleness-and-connectedness}
the map $\pi_1(Y_1 \cap Y_2) \to \pi_1(Y_1) \times \pi_1(Y_2)$
is surjective provided
$$
\dim Y_1+\dim Y_2> \dim X+ \dim \grass(\p^{r-1}, \p^{n-2}).
\eqno{(*)}
$$
Here $\dim X$ is what we expect on the right hand side of ($*$) from the naive dimension count. 
For $r=0$, $\dim \grass(\p^{r-1}, \p^{n-2}) = 0,$ and
we recover the Fulton-Hansen  theorem. However,  for $r\gg 0$ the extra term gets quite large. Nonetheless, this bound is optimal, see Example~\ref{non-exmp.2}.

Our aim is to prove that, for many group actions,
 the naive  dimension estimate $\dim Y_1+\dim Y_2> \dim X $ is enough
to give $\pi_1$-surjectivity,  unless the intersection $Y_1\cap Y_2$ is very degenerate. 
 
\begin{notation} \label{not.1.not}
We work over an algebraically closed field $k$ unless otherwise specified. 
Let $G$ be a connected algebraic group acting  on a normal $k$-variety $X$.
We usually assume that the action is faithful.

The action is denoted by $(g,x)\mapsto gx$. 
If $p:Y\to X$ is a morphism then it is convenient to use the shorthand
$gY\to X$ for $g\circ p:Y\to X$.

We consider the case when the $G$-action has a dense open orbit, denoted by  $X^\circ$, and also
a fixed point $x_0\in X$.
For any $n \geq 0$, $G$ acts linearly on the $n$th order infinitesimal neighborhood
$\mathscr{O}_{X,x_0}/\mathfrak{m}_{x_0}^n$, where $(\mathscr{O}_{X,x_0},
\mathfrak{m}_{x_0})$ is the
local ring at $x_0$ and its maximal ideal.
If we assume the action of $G$ on $X$ is faithful,
this induced action is also faithful for sufficiently large
$n$, and so, in this case, $G$ is a linear algebraic group.

We use $\pi_1(Z, z_0)$ to denote the 
topological fundamental group of $Z(\c)$ if we are over $\c$ and the 
\'etale fundamental group of $Z$ otherwise.
We will usually include basepoints, but sometimes omit them when it is not
important to keep track of them.
\end{notation}

The main result is easier to state for normal subvarieties over $\c$.
We prove a more general version in Theorem~\ref{dense.orb.thm}, where the
characteristic is arbitrary, the stabilizer may be disconnected  and the 
$Y_i$ are not necessarily subvarieties.

\begin{thm} \label{dense.orb.intro.thm}
Let $G$ be a connected algebraic group acting on a normal variety $X$ defined
over $\c$. Assume that the action has  a fixed point $x_0\in X$ and a dense open orbit $X^\circ\subset X$ 
such that the stabilizer $H \subset G$ of any point in $X^\circ$ is connected.
Let $Y_1, Y_2$ be irreducible, normal, locally closed
subvarieties of $X$ containing $x_0$. Assume that
\begin{enumerate}
\item  $\dim Y_1+\dim Y_2>\dim X$ and
\item there is an irreducible component  $Z\subset Y_1\cap Y_2$ 
of the expected dimension $\dim Y_1+\dim Y_2-\dim X$,
that contains $x_0$ and is {\em not}  disjoint from $X^\circ$.
\end{enumerate}
Then,
$$
\pi_1(g_1Y_1\cap g_2Y_2, x_0)\onto \pi_1(Y_1, x_0)\times \pi_1(Y_2, x_0) \qtq{is surjective}
$$
for general $(g_1, g_2)\in G\times G$.
\end{thm}

Note that the conclusion of Theorem~\ref{dense.orb.intro.thm}  is pretty much
the strongest that one can hope for:
If  $Y_1\cap Y_2$ 
is disjoint from $X^\circ$, then 
$g_1Y_1\cap g_2Y_2$  is never  in  `general position,' and if an irreducible component does not contain the base point $x_0$, then its contribution to the fundamental group  is hard to control; see also Examples~\ref{non-exmp.1}--\ref{non-exmp.2}.
For counterexamples with disconnected stabilizers see
Examples~\ref{example:nontrivial-fundamental-group}--\ref{example:disconnected-stabilizers}.

We next state and prove Theorem~\ref{dense.orb.thm}, from which
Theorem~\ref{dense.orb.intro.thm} follows.
We then derive generalizations of the Fulton-Hansen theorem
and further applications.
We conclude with various examples and counterexamples.

\begin{ack}   We thank
Brian Conrad and Zhiwei Yun for explaining the proof of 
Proposition~\ref{proposition:affine-to-projective}, 
Robert Lazarsfeld for pointing out Corollary~\ref{corollary:fulton-hansen-diagonal},
Laurent Manivel,
Anand Patel, and Arpon Raksit for helpful conversations, and
Akihiko~Yukie for long  comments and references on prehomogeneous vector spaces.
We also thank
Fran\cedilla{c}ois Charles, Dougal Davis, 
Laurent Moret-Bailly,
Wenhao Ou, 
Eric Riedl, Zev Rosengarten, 
Jason Starr,
Ashvin Swaminathan,
and
Ravi Vakil
for helpful correspondence on an earlier version of this article.

This is a preprint of an article published in 
Bollettino dell'Unione Matematica Italiana. The final
authenticated version is available online at: 
\newline
\url{https://doi.org/10.1007/s40574-021-00302-8}.

JK  was supported by  the NSF under grant number
DMS-1901855 and
AL  by the National Science Foundation Graduate Research Fellowship Program
under Grant No. DGE-1656518 and also supported through the program "Oberwolfach Research Fellows" by the Mathematisches Forschungsinstitut Oberwolfach in 2021. 
\end{ack}

\subsection*{Proof of the main theorem}{\ }

We start with some preliminary lemmas.

\begin{lem}
	\label{connected.fibers.lem}
Let $f:U\to V$ be a dominant morphism of irreducible varieties and
$\sigma:V\to U$ a section. Then $U_v$ is connected for general $v\in V$.
\end{lem}
\medskip 
Proof.
It is enough to show that the geometric generic fiber of $f$ is connected.
By irreducibility of $U$, the generic fiber of $f$ is irreducible.
The generic fiber may fail to be geometrically irreducible,
but all geometric components must pass through the image of $\sigma$,
and hence the geometric generic fiber is connected.
\qed
\medskip

We also require the following slight variant of 
\cite[Theorem 2]{kleiman:the-transversality-of-a-general-translate},
which follows from \cite[Theorem
2]{kleiman:the-transversality-of-a-general-translate}
and upper semicontinuity of dimension of fibers.

\begin{lem}\label{dim.fib.prod.lem}
Let $G$ be a connected algebraic group acting transitively on a normal
$k$-variety $X$.  Let $Y_1$ and $Y_2$ be irreducible, normal varieties  and 
 $p_i:Y_i\to X$ quasi-finite morphisms. Then there is a dense open subset
$U\subset G\times G$ such 
for all $(g_1, g_2)\in U$,
$g_1Y_1\times_X g_2Y_2$  is
\begin{enumerate}
\item either empty
\item or of pure dimension $\dim Y_1+\dim Y_2-\dim X$
\end{enumerate}
Moreover,  if $g_1Y_1\times_X g_2Y_2$  is non-empty and of dimension $\leq \dim Y_1+\dim Y_2-\dim X$ for some $(g_1, g_2)\in G\times G$, then (2) holds. \qed
\end{lem}

\begin{notation} \label{good.limit.defn}
Let $X$ be a $k$-variety, $x_0\in X$ a closed point and
$X^\circ\subset X$ a dense open subset. Let $p:Z\to X$ be a quasi-finite morphism. 
We let $Z^\circ$ denote $p^{-1}(X^\circ)$.
For $z_0\in p^{-1}(x_0)$, let $Z(z_0)\subset Z$  denote the union of those irreducible  components that contain $z_0$ and whose images are {\em not}  disjoint from $X^\circ$. We say that $z_0$ is a  {\it good limit point} of $Z$ if
  $Z(z_0)\neq\emptyset$. 

\end{notation}

\begin{lem} \label{good.lims.lem}
Let $G$ be a connected algebraic group acting on a normal variety $X$ with a
fixed point $x_0\in X$ and a dense open orbit $X^\circ$ with  stabilizer
$H\subset G$.  Let $H_0\subset H$ denote the identity component of $H$.

Let $Y_1$ and $Y_2$ be irreducible, normal varieties  and 
 $p_i:Y_i\to X$ quasi-finite morphisms such that  $Y_i^\circ:=p_i^{-1}(X^\circ)$ are nonempty.

Let $I\subset G \times G\times Y_1\times Y_2$ parametrize quadruples
$(g_1, g_2, y'_1, y'_2)$ such that $g_1^{-1}p_1(y'_1)=g_2^{-1}p_2(y'_2)$.
 Let $I^\circ\subset I$ be the preimage of 
$Y_1^\circ \times Y_2^\circ$ and  $\bar I^\circ\subset I$ its closure. We have  
$$
G \times G \stackrel{\pi}{\leftarrow} \bar I^\circ\stackrel{\tau}{\to} Y_1\times Y_2.
\eqno{(\ref{good.lims.lem}.1)}
$$
Then
\begin{enumerate}\setcounter{enumi}{1}
\item  
$I^\circ\to Y_1^\circ \times Y_2^\circ$ is a locally trivial $G\times H$-bundle.
\item $I^\circ $ has at most  $|H/ H_0|$ irreducible components.
\item Every good limit point $z_0=(y_1, y_2)$ 
	lying on a component of $Y_1 \times_X Y_2$ of dimension $\dim Y_1 + \dim Y_2 - \dim X$	
	gives a rational section $\pi_{z_0}: (g_1, g_2)\mapsto  \bigl(g_1, g_2, y_1, y_2)$. 
	\end{enumerate}
Let $j$ index those irreducible components
$\bar I^\circ(z_0, j)\subset \bar I^\circ$ that contain the image of $\pi_{z_0}$.
Then the following hold.
\begin{enumerate}\setcounter{enumi}{4}
\item For any $j$ as above, the image of $\pi_1\bigl(\bar I^\circ(z_0, j)\bigr)\to \pi_1(Y_1 \times Y_2)$ has index $\leq |H/ H_0|$.
\item $\bar I^\circ(z_0, j)\cap \pi^{-1}(g_1, g_2)$ is a  nonempty union of
	 irreducible components of $(g_1Y_1\times_X g_2Y_2)(z_0)$ for general $(g_1, g_2)\in G\times G$.
\item If $H$ is connected then $\pi^{-1}(g_1, g_2)$ is connected for general $(g_1, g_2)\in G\times G$.
\end{enumerate}
\end{lem}

{\it Remark \ref{good.lims.lem}.8.} It can happen that $g_1Y_1\times_X g_2Y_2 $ is always reducible, though this seems to be rare; see Example~\ref{reducible.exmp}.
\medskip

Proof.  Claims (2) and (4) are clear and (2) implies (3). There is nothing further to prove if there are no good limit points.

Otherwise let  $z_0=(y_1, y_2)$ be a good limit point.  
Set $K=k(G \times G)$. 
We next prove (6).
By Lemma~\ref{dim.fib.prod.lem},
the generic fiber $\bar I^\circ(z_0, j)_K$ of 
the projection $\bar I^\circ(z_0, j)\to G \times G$ is nonempty, hence irreducible over $K$, and
it has a $K$-point, given by $z_0$.  
Note there is some such component $\bar I^\circ(z_0, j)$ because $z_0$ is a good
limit point.
The union of those geometric irreducible components of $\bar I^\circ(z_0, j)_K$ that contain this $K$-point is defined over 
an extension of $K$.
Thus, every geometric irreducible component of $\bar I^\circ(z_0, j)_K$ contains 
this  $K$-point. Hence,  for general $(g_1, g_2)\in G\times G$, 
every  irreducible component of  $\bar I^\circ(z_0, j)\cap\pi^{-1}(g_1, g_2)$
contains  $(g_1, g_2, z_0)$, proving (6).

As for (5),  $I^\circ\to Y_1^\circ \times Y_2^\circ$ factors as
$$
I^\circ\stackrel{\tau_1}{\to} I^\circ/(G\times H_0)\stackrel{\tau_2}{\to} Y_1^\circ \times Y_2^\circ.
$$
Here $\tau_1$ is a principal $G \times H_0$-bundle with connected fibers and
$\tau_2$ is finite, \'etale of degree $|H/ H_0|$ (but possibly disconnected).
Thus the image of
$$
\pi_1(I^\circ)\to \pi_1(Y_1^\circ \times Y_2^\circ)
$$
is a subgroup of index $\leq |H/ H_0|$. 
Observe that
$\pi_1(Y_1^\circ \times Y_2^\circ)\to \pi_1(Y_1 \times Y_2)$ is surjective 
as $Y_1$ and $Y_2$ are
normal; see \cite[0.7.B]{fultonL:connectivity-and-its-applications}.  
We can thus factor $\tau: \bar I^\circ(z_0, j)\to Y_1 \times Y_2$ as
$$
\bar I^\circ(z_0, j)\xrightarrow{\tau_0}(Y_1 \times Y_2)^{\sim}
\xrightarrow{\tilde\tau} Y_1 \times Y_2,
\eqno{(\ref{good.lims.lem}.9)}
$$
where $\tilde\tau$ is finite, \'etale of degree $\leq |H/ H_0|$
and
$$
\pi_1\bigl(\bar I^\circ(z_0, j)\bigr)\to \pi_1\bigl((Y_1 \times Y_2)^\sim\bigr)\qtq{is surective.}
\eqno{(\ref{good.lims.lem}.10)}
$$
This proves (5) and  (7) follows from Lemma~\ref{connected.fibers.lem}.
\qed

\medskip

The proof of Theorem~\ref{dense.orb.thm} will follow from the above, together
with the following important fact from topology.
\begin{rem}
	\label{remark:fiber-bundle}
	Suppose we are given a map $f: X \to Y$ of varieties over $\c$.
	We claim that there is a Zariski open subset $U \subset Y$ such that 
	$X \times_Y U \to U$ is locally a trivial bundle in the Euclidean
	topology. 
	More precisely, 
	for any point $p \in U^{\an}$, the analytic space associated to $U$,
	there is an analytic open set $V \subset U^{\an}$ with $p \in V$
	and a homeomorphism of topological spaces
	spaces $(f^{\an})^{-1}(V) \simeq (f^{\an})^{-1}(p) \times V$ over $V$.

	The claim above follows from
	\cite[Part I, Theorem 1.7, p. 43]{goresky1988stratified}
	and its proof, although it is not explicitly stated there,
	so we comment on the details.
First, \cite[\S4, Theorem 1]{hironaka:stratification-and-flatness}
shows that the pieces in the stratification of
\cite[Part I, Theorem 1.7, p.\ 43]{goresky1988stratified}
can be taken to be complex analytic. 
(Note that the pieces of the stratification are not necessarily real analytic when the map is real analytic, 
see \cite[Caution, p.\ 43]{goresky1988stratified}.)
To prove the existence of the desired Zariski open set $U \subset Y$, one may work with proper
compactifications $\overline{X}, \overline{Y}$ of $X$ and $Y$,
which exist by Nagata's compactification theorem.
We may then use the last sentence of
\cite[Part I, Theorem 1.7, p.\ 43]{goresky1988stratified}
to ensure the boundaries of these closures are also unions of strata of the stratifications. 
Hence, we also obtain stratifications of the open parts $X$ and $Y$.
Then, \cite[Part I, \S1.6, p.\ 43]{goresky1988stratified}
explains why the resulting map is topologically trivial over pieces of the
stratification.
Finally, to show that one of the pieces of the stratification of $Y$ is Zariski open,
we note that there is a piece whose complement is a complex analytic proper closed subvariety, hence algebraic by Chow's theorem.
\end{rem}

\begin{say}[Fundamental groups of fibers] \label{pi1.fiber.say}
Let $S, W, Z$ be irreducible $k$-varieties, $Z$ is normal. Assume that we have
 dominant morphisms
$$
S\stackrel{q}{\longleftarrow} W \stackrel{p}{\longrightarrow} Z
$$
and  a section  $\sigma:S\to W$ of $q$ such that $p\circ \sigma:S\to Z$ is a constant morphism with image $z_0\in Z$. 
Assume that $\pi_1(W)\to \pi_1(Z)$ is surjective and $W \to Z$ is a flat map
with irreducible fibers.  We aim to understand whether the induced homomorphism
$\pi_1\bigl(W_s, \sigma(s)\bigr)\to \pi_1(Z, z_0)$ is surjective for general  $s\in S$.
\medskip

{\it Case \ref{pi1.fiber.say}.1.}  $k=\c$. By Remark~\ref{remark:fiber-bundle}, there is an open, dense  subset $S^\circ\subset S$ such that   $ q$ is a topological  fiber bundle over $S^\circ$. Set $W^\circ:=q^{-1}(S^\circ)$. Note that 
$W^\circ\to S^\circ$   has a section that gets contracted by $p$.  Thus
$$
\pi_1\bigl(W_s, \sigma(s)\bigr)\onto \pi_1(Z, z_0) \qtq{is surjective for every $s\in S^\circ$.}
\eqno{(\ref{pi1.fiber.say}.1.a)}
$$

\medskip

{\it Case \ref{pi1.fiber.say}.2.}  $\chr k=0$. By the Lefschetz principle we get
that (\ref{pi1.fiber.say}.1.a) holds for the \'etale fundamental group. 

\medskip

{\it Case \ref{pi1.fiber.say}.3.}  $\chr k>0$.  Let  $\pi_Z: (Z', z'_0)\to (Z, z_0)$ be an irreducible,  finite, \'etale cover  corresponding to a finite quotient
$\pi_1(Z, z_0)\onto H$. By pullback we obtain a
connected,  finite, \'etale cover  $\pi_W:W'\to W$. 
Here, $W'$ is  irreducible
because $W' \to Z'$ is flat with irreducible fibers.
Further, $\sigma$ lifts (non-uniquely) to a section $\sigma':S\to W'$ such that $p'\circ \sigma'$ maps $S$ to $z'_0$. 
Thus, by Lemma \ref{connected.fibers.lem}, there is an open, dense  subset $S^\circ_H\subset S$ such that  $W'_s$ is connected for $s\in S^\circ_H$.  
That is,
$$
\pi_1\bigl(W_s, \sigma(s)\bigr)\to \pi_1(Z, z_0)\to H \qtq{is surjective  for every $s\in S^\circ_H$.}
\eqno{(\ref{pi1.fiber.say}.3.a)}
$$
In this case we say that  $\pi_1\bigl(W_s, \sigma(s)\bigr)\to \pi_1(Z, z_0)$ is  {\it surjective on finite quotients} for general $s\in S$.

\medskip
{\it Remark \ref{pi1.fiber.say}.4} The complication in
\ref{pi1.fiber.say}.3
is that $S^\circ_H$ depends on $H$, and, if $\chr k>0$,then 
$\cap_H S^\circ_H$ may be empty. That is,
there may not be any closed point $s\in S$ for which
$\pi_1\bigl(W_s, \sigma(s)\bigr)\to \pi_1(Z, z_0)$  is surjective; see
Example~\ref{example:positive-characteristic}. 
\end{say}

Combining Lemma \ref{good.lims.lem} with Paragraph \ref{pi1.fiber.say}, 
can now prove our main theorem.

\begin{thm} \label{dense.orb.thm}
Let $G$ be a connected algebraic group acting on a normal $k$-variety $X$ with a
fixed point $x_0\in X$ and a dense open orbit $X^\circ$  with  stabilizer
$H\subset G$.  Let $(y_1,Y_1), (y_2,Y_2)$ be irreducible, normal, pointed  varieties  and 
 $p_i:Y_i\to X$ quasi-finite morphisms such that $p_i(y_i)=x_0$. 
Assume that 
\begin{enumerate}
\item  $1\leq \dim (Y_1^\circ\times_X Y_2^\circ)\leq \dim Y_1+\dim Y_2-\dim X$, and 
\item $Y_1\times_X Y_2$ has a good limit point
	$z_0=(y_1,y_2)$, as in Notation~\ref{good.limit.defn}.
 \end{enumerate}
Then,  for general $(g_1, g_2)\in G\times G$,  
\begin{enumerate}\setcounter{enumi}{2}
\item  the closure of $g_1Y_1^\circ\times_X g_2Y_2^\circ$ in $g_1Y_1\times_X g_2Y_2$ has at most
$|H/ H_0|$ connected components, and 
\item  there is a subgroup  $\Gamma\subset\pi_1(Y_1, y_1)\times \pi_1(Y_2, y_2)
	$ of index at most $|H/ H_0|$ such that the natural map, 
$$
\pi_1\bigl(g_1Y_1\times_X g_2Y_2, z_0\bigr)\to \Gamma
$$
is  surjective if $\chr k=0$, and  surjective  on finite quotients if $\chr k>0$.
\end{enumerate}
\end{thm}
Proof.   Fix an index $j$ as in  Lemma~\ref{good.lims.lem} and
consider the diagram
$$
G \times G \xleftarrow{\pi} \bar I^\circ(z_0, j)\xrightarrow{\tau_0}(Y_1 \times Y_2)^\sim
\xrightarrow{\tilde\tau} Y_1 \times Y_2,
$$
whose right hand side is defined in (\ref{good.lims.lem}.9).
Set 
$$
\Gamma:=\pi_1\bigl((Y_1 \times Y_2)^\sim\bigr)\subset\pi_1(Y_1, y_1)\times \pi_1(Y_2, y_2).
$$
Next apply the discussions in Paragraph~\ref{pi1.fiber.say} to conclude that
$$
\pi_1\bigl(\bar I^\circ(z_0, j)\cap \pi^{-1}(g_1, g_2)\bigr)\to \pi_1\bigl((Y_1 \times Y_2)^\sim\bigr)
$$
is  surjective if $\chr k=0$, and  surjective  on finite quotients if $\chr k>0$. \qed

\subsection*{Fulton-Hansen-type  theorems}{\ }

We get the following version of the Fulton-Hansen theorem.

\begin{cor}[Fulton-Hansen theorem I] \label{corollary:fulton-hansen-diagonal}
Let $\mathcal G$ be a connected algebraic group acting 2-transitively 
on a  quasi-projective variety $X$.
Let  $Z$ be a normal, irreducible variety with $\dim Z > \dim X$ and $p: Z \to X \times X$
quasi-finite. 
Then,  for general $(g_1, g_2)\in \mathcal G\times \mathcal G$, and $\Delta_X :X
\to X \times X$ the diagonal map,
\begin{enumerate}
	\item  $Z \times_{(g_1, g_2) \circ p, X\times X, \Delta_X} X$ is connected and 
\item  the natural map, 
$$
\pi_1\bigl(Z \times_{(g_1, g_2) \circ p, X \times X, \Delta_X} X, z_0\bigr)\onto \pi_1(Z,z)\times
\pi_1(X,x),
$$
is  surjective if $\chr k=0$ and surjective on finite quotients if $\chr k>0$.
\end{enumerate}
\end{cor}

Proof (assuming Corollary~\ref{lemma:2-orbit-classification.cor}).
For general $(g_1, g_2)\in \mathcal G\times \mathcal G$ 
choose 
$z_0 = (z,x) \in Z \times_{(g_1, g_2) \circ p, X \times X, \Delta_X} X$ and define $x_0
\in X$
as the point such that
$(x_0, x_0)=((g_1, g_2)\circ p)(z)$.
Let $G\subset \mathcal G$ denote the identity component of the stabilizer of
$x_0$. Note that $G$ acts transitively  on $X\setminus\{x_0\}$.
(By Remark~\ref{2-trans.actions.rem}, the stabilizer in $\mathcal G$ of $x_0$ is
already connected, though we will not need this fact.)

We aim to apply
Theorem~\ref{dense.orb.thm} to  $G\times G$ acting on $X \times X$ with 
fixed point $(x_0, x_0)$ and dense
open orbit $(X\setminus \{x_0\}) \times (X \setminus \{x_0\})$.
We now check the hypotheses.
As an itinitial reduction, the statement for general translations of both $Z$ and $X$ is equivalent to 
the analogous one for only translates of $Z$ because
$Z \times_{(g_1, g_2) \circ p, X \times X, (h_1, h_2) \circ \Delta_X} X
\simeq Z \times_{(h_1^{-1} g_1, h_2^{-1} g_2) \circ p, X \times X,
\Delta_X} X$.
Now, note that, $\dim_{z_0}(Z \times_{(g_1, g_2) \circ p, X \times X, \Delta_X} X)=\dim Z + \dim X
- \dim X \times X> 0$ by
Lemma~\ref{dim.fib.prod.lem}.
Since $\Delta_X(X)\setminus\{(x_0,x_0)\}$ is contained in the dense orbit $(X
\setminus\{x_0\}) \times (X \setminus \{x_0\})$,
$z_0$ is a good limit point of $Z \times_{(g_1, g_2) \circ p, X \times X, \Delta_X} X$.
Additionally, it follows from Corollary~\ref{lemma:2-orbit-classification.cor} that 
the stabilizer
$H \subset G$ of  the $G$ action on $X\setminus\{x_0\}$  is connected. 
Therefore, the stabilizer in $G \times G$ of a point of 
$X\setminus\{x_0\} \times X\setminus\{x_0\}$ 
is also connected.
\qed

\medskip

We could have proven Corollary~\ref{corollary:fulton-hansen-diagonal}
using Remark~\ref{2-trans.actions.rem} in place of 
Corollary~\ref{lemma:2-orbit-classification.cor}.
We opted to use the latter as it leads to 
 a more self-contained proof.

An important special case of Corollary~\ref{corollary:fulton-hansen-diagonal}
is the the case where $Z = Y_1 \times Y_2$ and $p = p_1 \times p_2$ for $p_i:
Y_i \to X$ quasi-finite morphisms.
Because this is typically how Corollary~\ref{corollary:fulton-hansen-diagonal} is
applied, we now restate it in this case.

\begin{cor}[Fulton-Hansen theorem II] \label{ful-han.thm.cor}
Let $\mathcal G$ be a connected algebraic group acting 2-transitively 
on a  quasi-projective variety $X$.
Let  $Y_1, Y_2$ be normal, irreducible varieties and $p_i:Y_i\to X$ quasi-finite morphisms. Assume that $\dim Y_1+\dim Y_2>\dim X$. 
Then,  for general $(g_1, g_2)\in \mathcal G\times \mathcal G$,  
\begin{enumerate}
\item  $g_1 Y_1\times_X g_2 Y_2$ is connected and 
\item  the natural map, 
$$
\pi_1\bigl(g_1Y_1\times_X g_2Y_2, z_0\bigr)\onto \pi_1(Y_1, y_1)\times
\pi_1(Y_2, y_2),
$$
is  surjective if $\chr k=0$ and  surjective  on finite quotients if $\chr k>0$.
\end{enumerate}
\end{cor}
%Proof (assuming Corollary~\ref{lemma:2-orbit-classification.cor}).
%For general $(g_1, g_2)\in \mathcal G\times \mathcal G$ 
% there is  a point
%$z_0=(y_1, y_2)\in g_1 Y_1 \times_X g_2 Y_2$ such that $Y_i$  are normal at
%$y_i$. Set $x_0=p_1(g_1 y_1)=p_2(g_2 y_2)$ and
%let $G\subset \mathcal G$ denote the identity component of the stabilizer of
%$x_0$. Note that $G$ acts transitively  on $X\setminus\{x_0\}$.
%(By Remark~\ref{2-trans.actions.rem} the stabilizer $G$ is 
%connected, though we will not need this fact.)
%
%We aim to apply
%Theorem~\ref{dense.orb.thm} to  $G$ acting on $X$,
%and we now check the hypotheses.
%First, $\dim_{z_0}(g_1Y_1\times_Xg_2Y_2)=\dim Y_1+\dim Y_2-\dim X$ by
%Lemma~\ref{dim.fib.prod.lem}.
%Thus, $z_0$ is a good limit point of $g_1Y_1\cap g_2Y_2$.
%Finally, it follows from Corollary~\ref{lemma:2-orbit-classification.cor} that 
%the stabilizer
%$H \subset G$ of  the $G$ action on $X\setminus\{x_0\}$  is connected. 
%(We could have proven Corollary~\ref{ful-han.thm.cor}
%using Remark~\ref{2-trans.actions.rem} in place of 
%Corollary~\ref{lemma:2-orbit-classification.cor}; the latter leads to 
% a more self-contained proof.)\qed
%
\medskip

If $X$ and $Z$ as in Corollary~\ref{corollary:fulton-hansen-diagonal} are proper, then connectedness of 
$Z \times_{(g_1, g_2) \circ p, X \times X, \Delta_X} X$  for general $(g_1, g_2)\in \mathcal G\times \mathcal G$ implies connectedness for every $(g_1, g_2)$. This is sometimes called the Enriques-Severi-Zariski connectedness principle, proved by combining Stein factorization with Zariski's main theorem.

Thus we recover the original setting of the 
Fulton-Hansen therem \cite{fultonH:a-connectedness-theorem}. 

\begin{cor}[Fulton-Hansen theorem III] \label{ful-han.thm.cor-ii}
Let $\mathcal  G$ be a connected algebraic group acting 2-transitively on a projective
variety $X$. 
Let $Z$ be a normal irreducible proper variety with $\dim Z > \dim X$ and $p: Z \to X \times X$ a
finite morphism.
Then $Z \times_{p, X\times X, \Delta_X} X$ is connected.
In particular, if $Z = Y_1 \times Y_2$ and $p = p_1 \times p_2$ for $p_i: Y_i
\to X$ finite morphisms, then
$Y_1\times_{p_1, X, p_2} Y_2$ is connected. \qed
\end{cor}

\begin{rem}\label{2-trans.actions.rem}
 The most important example of a 2-transitive action is
$(\PGL_{n+1}, \p^n)$, 

The 2-transitive case seems much more general, but in fact there are very few
such pairs $(\mathcal G, X)$. By \cite{knop:mehrfach-transitive}
$(\PGL_{n+1}, \p^n)$ is the only pair with $X$ projective. 
The pairs with $X$ quasi-projective are
all of the form  $(G\ltimes \ga^n, \a^n)$ where $G\subset \GL_n$ is a product
$C \cdot \overline{G}$, 
for
$C$ a subgroup of the central $\mathbb G_m \subset \GL_n$ and
$\overline{G}$ is one of the
following:
\begin{enumerate}
\item $n = 1$, $\overline{G} = \GL_1$, 
\item $n \geq 2$, $\overline{G} = \SL_n$,
\item $n=2m$ is even, $\overline{G} = \Sp_{m}$,
\item $n =6$, the characteristic is $2$, and  $\overline{G} = \operatorname{G}_2$. 
\end{enumerate}
(Note that $\operatorname{G}_2 $ does not have a nontrivial 6-dimensional representation in  characteristics  $\neq 2$.)
\end{rem}

There is, however, a very long list of pairs  $(G, X)$ such that
$G(\r)$ acts 2-transitively on $X(\r)$; see
\cite{tits:sur-certaines-classes, kramer:two-transitive-lie-groups}.
So the following variant applies in many more cases.

\begin{cor} Let $X$ be a variety defined over $\r$ and $G$  a connected algebraic group acting on it
such that the $G(\r)$ action on   $X(\r)$ is 2-transitive. 
Assume that for $x_0\neq x_1\in X(\r)$, the  stabilizer of the ordered pair $(x_0, x_1)$ is connected (over $\c$). 
Let  $Z_1, Z_2$ be irreducible, normal  varieties and $p_i:Z_i\to X$ quasi-finite morphisms.
Assume that $\dim Z_1+\dim Z_2>\dim X$ and the $Z_i$ have smooth real points.
Then,
$$
\pi_1(g_1Z_1\times_X g_2Z_2)\to  \pi_1(Z_1)\times\pi_1(Z_2)\qtq{is surjective}
$$
for general $(g_1, g_2)\in G(\r) \times G(\r)$ and  for general $(g_1, g_2)\in G(\c)
\times G(\c)$. \qed
\end{cor}

There are also some non-transitive group actions for which 
we get a Fulton-Hansen-type  result, with obvious exceptions.

\begin{exmp}[Orthogonal group]\label{orth.exmp} Let $\GO_q:=\gm\cdot \OO_q$ be the 
group of orthogonal similitudes acting on the $n$-dimensional vector space $V^n$, where $q$ is a nondegenerate quadratic form. 
There are 3 orbits, $\{0\}$,  $(q=0)\setminus\{0\}$, and
the dense open orbit  is
$V^n\setminus (q=0)$. 
\medskip

{\it Claim \ref{orth.exmp}.1.} Let $0\in Y_i\subset \a^n$ be irreducible,
normal, locally closed subvarieties. Assume that $Y_i\not\subset (q=0)$ and 
$\dim Y_1+\dim Y_2>n$.  
Then 
$$
\pi_1(g_1Y_1\cap g_2Y_2,0)\onto \pi_1(Y_1,0)\times \pi_1(Y_2,0) \qtq{is surjective}
$$
for general $(g_1, g_2)\in  \GO_q\times\GO_q$.

\medskip

Proof. Since $\a^n$ is smooth, 
$$
\dim_0 (g_1Y_1\cap g_2Y_2)\geq \dim Y_1+\dim Y_2-n.
$$
Thus we have a good limit point if
$$
\dim_0 \bigl(g_1Y_1\cap g_2Y_2\cap (q=0)\bigr)< \dim Y_1+\dim Y_2-n.
$$
 Since $Y_i\not\subset (q=0)$,  we see that
$\dim_0 \bigl(g_iY_i\cap (q=0)\bigr)\leq  \dim Y_i-1$. 
Since $(q=0)\setminus\{0\}$ is homogeneous, using Lemma~\ref{dim.fib.prod.lem} we see that
$$
\dim_0 \bigl(g_1Y_1\cap g_2Y_2\cap (q=0)\bigr)\leq 
(\dim Y_1-1)+(\dim Y_2-1)-(n-1),
$$
as needed. \qed
\end{exmp}

The above arguments show that our approach gives  the best results
if the orbits of an action are fully understood. In the most extreme case, we have the following classification. The proof relies on some results of
\cite{knop:mehrfach-transitive}, that we recall afterwards.

\begin{prop}
	\label{lemma:2-orbit-classification}
	Let $X$ be an irreducible,  normal variety of dimension $\geq 2$ over a   field $k$	and  $G$
	 a  connected linear algebraic group  acting on $X$.
Assume that all orbits have dimension either $0$ or $\dim X$.
\begin{enumerate}
\item There is at most 1 orbit of dimension $0$.
\item If  $\chr k = 0$ and there is a $0$-dimensional orbit, then
  $X$ is isomorphic to either an affine or a projective cone
over a projective, homogeneous $G$-variety  $Y$.
\end{enumerate}
\end{prop}

Proof.  
For (1) we may assume that $k$ is algebraically closed.
We may then replace $G$ by its reduction to assume $G$ is smooth.
Let $P=\{p_i\}$ be the union of the 0-dimensional orbits and assume $P$ is
nonempty.
By Proposition~\ref{proposition:affine-to-projective}
there is a projective $G$-variety $Y$ and a $G$-equivariant,  affine, surjective morphism $f: X\setminus P \to Y$, whose general fiber is 1-dimensional 
by Lemma~\ref{lemma:relative-dimension-1}.

By \cite[Theorem 3]{sumihiro:equivariant-completion}, there is a normal, 
$G$-equivariant compactification  $\bar X\supset  X$. 
Let $Z$ be the normalization of the closure of the graph of $f$ with projections  $\pi_X$ and $\pi_Y$. 

Since $G$ acts transitively on $X\setminus P$, it also acts
transitively on $Y$, hence $E:=Z\setminus   (X\setminus P)$
is a union of $G$-orbits.
Thus every fiber of $\pi_Y:Z\to Y$ is  a geometrically rational curve
 and $E$ is a disjoint union of (possibly multiple) sections.

For any $p_i\in P$ there is an irreducible component  $E_i\subset E$ that is contracted by $\pi_X$ to $p_i$. Let  $C\subset Y$ be a general curve, 
$X_C$ the normalization of  $\pi_Y^{-1}(C)$   and $F_i\subset X_C$ the preimage of $E_i$.  Then $X_C\to C$ is a   $\p^1$-bundle. Note that  a $\p^1$-bundle over a smooth, projective curve contains at most 1 curve with negative self-intersection, and this curve is a section. Thus $P$ has at most 1 point.

In order to prove (2), we assume from now on that $\chr k=0$. 
Then $\pi_Y:Z\to Y$ is  a $\p^1$-bundle and 
 $E$ consists of 1 or  2 sections.

If $E$ consists of 2 sections, then $Z=\p_Y(\o_Y+L)$ for some anti-ample line bundle  $L$  on $Y$. Thus $X$ is the affine cone  over $(Y, L^{-1})$. 

If $E$ consists of 1 section, then $Z=\p_Y(E)$  where $E$ is obtained as an extension
$$
0\to \o_Y\to E\to L\to 0,
$$
for some anti-ample line bundle  $L$  on $Y$.
Now we again use that  $\chr k=0$, hence $-K_Y$ is ample
and Kodaira's vanishing theorem implies that the extension splits. 
Thus $X$ is the projective cone  over $(Y, L^{-1})$. \qed
\medskip

{\it Remark \ref{lemma:2-orbit-classification}.3.}  We believe that (2) is not true if $\chr k\neq 0$, but  it should be possible to get a complete description of all cases.

\medskip
From this we deduce a useful corollary, used in the proof of
Corollary~\ref{ful-han.thm.cor} above.

\begin{cor}
	\label{lemma:2-orbit-classification.cor}
	Let $X$ be an irreducible,  smooth variety   over a   field $k$
	and  $G$
	 a  connected linear algebraic group  acting on $X$  with 2 orbits, one of which  is a	point $x \in X$.  
Then  $X\setminus\{x\}\cong G/H$ where $H$ is connected.
\end{cor}

Proof. This is clear if $\dim X=1$, so assume that $\dim X\geq 2$.

Let $H_0\subset H$ be the identity component.
If $H_0\neq H$ then $G/H_0\to G/H$ is an \'etale cover, which extends to an
\'etale cover $\pi: \tilde X\to X$ by purity. The $G$-action on $\tilde X$ has an open orbit  $\tilde X\setminus \pi^{-1}(x)$ and $\pi^{-1}(x)$ is a union of 0-dimensional orbits. Thus $\deg \pi=1$ by Propositon~\ref{lemma:2-orbit-classification}.  \qed

\begin{lem}
	\label{lemma:relative-dimension-1} {\rm (cf.\ \cite[p.\ 443]{knop:mehrfach-transitive})}
	Let $X$ be a normal irreducible variety, $P\subset X$ a 0-dimensional subset  and $f:X\setminus P \to Y$ an
	affine morphism.
	Then $\dim X\leq \dim Y+1$.
\end{lem}

Proof. Assume that $\dim X\geq \dim Y+1$. 
Choose normal compactifications  $\bar X\supset X$, $\bar Y\supset Y$
and let $Z$ be the normalization of the closure of the graph of $f$ with projections  $\pi_X$ and $\pi_Y$.   Note that $\pi_X$ cannot contract a whole fiber of $\pi_Y$. Thus there is point $\bar y\in \bar Y$ and an 
irreducible component  $Z_z\subset \pi_Y^{-1}(\bar y)$ such that
$Z_z\cap \pi_X^{-1}(P)$ and  $Z_z\setminus \pi_X^{-1}(P)$ are both nonempty. Set
$W=\pi_X(Z_z)$. Note that $P\cap W\neq\emptyset$, so $W^\circ:=W\cap X$ is dense in $W$.
Thus $ W^\circ\cap (X\setminus P)$ is the fiber of $f$, hence affine. 
By Hartogs's theorem, this implies $\dim W=1$. So the general fiber dimension of $f$ is $\leq 1$. \qed

\medskip

The following group theoretic result is proved, but not stated, on \cite[p.\
443]{knop:mehrfach-transitive}. 
It does not seem to be well known, so we now state it and give a proof in
Paragraph~\ref{subsection:knop-proof}, following suggestions of
Brian Conrad and Zhiwei Yun.

\begin{prop}
	\label{proposition:affine-to-projective}
	Suppose $X = G/H$ is a homogeneous space for a smooth connected linear algebraic group $G$
	over an algebraically closed field $k$.
	Then, there exists a parabolic subgroup $H \subset P \subset G$ with
	$P/H$ affine.
	In particular, there exists a projective variety $Y$ and a surjective and
	affine map
	$X \to Y$.
\end{prop}

\begin{rem} 
	\label{remark:canonical-choice-of-parabolic}	
	Knop's proof of Proposition~\ref{proposition:affine-to-projective} 
	in \cite{knop:mehrfach-transitive}
	is
	slightly different from ours in that he produces a specific choice of
	$P_H$ associated to $H$, whereas our proof merely takes $P$ to be an arbitrary minimal
	parabolic containing $H$.
	It is not clear to us what internal property distinguishes it from the other choices.

	An interesting aspect is that $P$ is usually not unique and the set of such parabolics has neither a smallest nor a largest element.
For example, for
$$
H:=\begin{pmatrix}
			* &  0   \\
			0 & * 
		\end{pmatrix}
\subset \SL_2
$$
the maximal choice of $P$ is $\SL_2$; the minimal choices are either the upper or the lower triangular matrices.
For 
$$
H:=\begin{pmatrix}
			1 &  0 & *  \\
			0 & 1 & 0 \\
0 & 0 & 1  		
		\end{pmatrix}
\subset \SL_3
$$
the minimal choice is $B:=(\mbox{upper triangular matrices})$. 
The maximal choices are the 2 maximal proper parabolics that contain $B$. 

\end{rem}

In what follows, for $G$ an algebraic group, we use $R_u(G)$ to denote its
unipotent radical, the maximal smooth normal connected unipotent subgroup of
$G$.
\begin{lem}
	\label{lemma:maximal-subgroup-reductive-or-parabolic}
	Let $H$ be a subgroup of a smooth connected reductive group $G$
	over a perfect field $k$.
	Then $H$ is either reductive or contained in a $k$-parabolic subgroup of
	$G$.
\end{lem}
Proof.
The key input in this proof is the fact that any smooth connected unipotent subgroup $U$ 
of a connected linear algebraic group $G$ over a perfect field $k$ 
is contained in the unipotent radical
of a parabolic $k$-subgroup $P \subset G$.
This follows from a theorem of Bruhat-Tits, see the ``Refined Theorem'' in
\cite{conrad:mathoverflow-refined-borel-tits}.

Applying this to our situation, suppose $H$ is not reductive.
We wish to show $H$
is contained in a proper parabolic subgroup.
By the above fact, there is some $P \subsetneq G$ with $R_u(H) \subset R_u(P)$.
It follows that $R_u(H) \subset H \subset N_G(R_u(H)) \subset
N_G(R_u(P)) \subset N_G(P) = P$.
Therefore, $H \subset P \subsetneq G$ for $P$ parabolic.
\qed
\medskip

\begin{say}[Proof of Knop's Proposition~\ref{proposition:affine-to-projective}.] 
	\label{subsection:knop-proof}
We can write $X = G/H$ for $H \subset G$ a subgroup.
Let $P$ denote a minimal parabolic containing $H$.
We wish to show $P/H$ is affine.
Let $L := P/R_u(P)$
and let $K$ denote the image of $H$ in $L$.
Because $P$ was chosen to be a minimal parabolic containing $H$, $K$ is not contained in any proper
parabolic subgroup of $L$.
By Lemma \ref{lemma:maximal-subgroup-reductive-or-parabolic},
$K$ is reductive.
Let $U := \ker(H \to K)$.
Then, $U \subset R_u(P)$ so $U$ is unipotent.
The quotient $R_u(P)/U$ is an affine group scheme acting on
$P/H$ with quotient $L/K$.
It follows that $P/H$ is a principal $R_u(P)/U$-bundle over $L/K$, and so to show
$P/H$ is affine, it suffices to show $L/K$ is.
This follows from the general claim that a quotient of a connected
reductive group by a connected reductive subgroup is affine, see
\cite[Theorem 1.5]{borel:on-affine-algebraic-homogeneous-spaces}.
\qed
\end{say}
\medskip

\subsection*{Other applications}{\ }

\begin{exmp}[Projective homogeneous spaces]\label{proj.hom.exmp}   These are of
	the form  $X=G/P$ where $G$ is a semisimple algebraic group and
	$P\subset G$ a parabolic subgroup. We get a Schubert cell decomposition with a single fixed point
$x_0$ and an open cell  $X^*\subset X$. Note that $X^*$ is a homogeneous space under the unipotent radical $U\subset P$.
The stabilizer of the $U$-action on $X^*$ is trivial, hence connected. Thus we get the following.
\medskip

{\it Claim \ref{proj.hom.exmp}.1.}  Let $Y_1, Y_2$ be irreducible, normal varieties  and  $p_i:Y_i\to X$ quasi-finite morphisms.
Assume that there is an irreducible component 
$$
Z^*\subset Y_1\times_X Y_2\times_X X^*\qtq{and a point} z_0=(y_1, y_2)\in Z\subset Y_1\times_X Y_2,
$$ 
such that $\dim Z^*=\dim Y_1+\dim Y_2-\dim X$ and  $p_i(y_i)=x_0$, where $Z$ denotes the closure of $ Z^*$.
Then, for general $(g_1, g_2) \in U \times U$, (hece also for general $(g_1, g_2) \in G \times G$),  the natural map
$$
\pi_1\bigl(g_1Y_1\times_X g_2Y_2, z_0\bigr)\to \pi_1(Y_1, y_1)\times \pi_1(Y_2, y_2) \qtq{is}
$$
\begin{enumerate}%\setcounter{enumi}{2}
\item[(a)] surjective if $\chr k=0$ and 
\item[(b)] surjective  on finite quotients if $\chr k>0$.
\end{enumerate}
\medskip

Note that we could consider instead the $P$-action, which has a
usually  larger open orbit  $X^\circ\supset X^*$. This gives the following variant.
\medskip

{\it Claim \ref{proj.hom.exmp}.2.}  Using the above notation, assume that
 $$
\dim Y_1+\dim Y_2> \dim X+ \dim (X\setminus X^\circ).
\eqno{(\ref{proj.hom.exmp}.2.a)}
$$
Then  the natural map
$$
\pi_1\bigl(g_1Y_1\times_X g_2Y_2, z_0\bigr)\onto \pi_1(Y_1, y_1)\times \pi_1(Y_2, y_2) \qtq{is}
$$
\begin{enumerate}%\setcounter{enumi}{2}
\item[(b)] surjective if $\chr k=0$, and 
\item[(c)] surjective  on finite quotients if $\chr k>0$,
\end{enumerate}
 for general $(g_1, g_2)\in P\times P$. 
\medskip

Proof. For dimension reasons  there is an irreducible component 
$Z^\circ\subset Y_1\times_X Y_2\times_X X^\circ$ that contains a good limit point $z_0$. Next we use the $P$-action to see that
$z_0$ is also a good limit point in  $g Z^\circ\times_X X^*$ for general $g\in P$.
 Now we can apply (1). \qed
\medskip

Results of this type have been considered in
\cite{sommese:complex-subspaces-i, faltings:formale-geometrie-und-homogene-raume,  goldstein:ampleness-and-connectedness, hansen:a-connectedness-theorem-for-flagmanifolds-and-grassmannians}.
Our bound (\ref{proj.hom.exmp}.2.a) is optimal in some cases, but  weaker in several of them. 
Using the full $G$-action, as in the above articles, leads to further improvements, but we did not find a natural way to recover the  bounds of \cite{faltings:formale-geometrie-und-homogene-raume, goldstein:ampleness-and-connectedness} in all cases.

\end{exmp}

\begin{exmp}[Prehomogeneous vector spaces]  
  A {\it  prehomogeneous vector space} is a pair  
$(G, V^n)$ where $V^n$ is a $k$-vector space of dimension $n$ and 
$G\subset \GL_n$ is a   connected subgroup that has a dense orbit  $W\subset
V^n$.   See \cite{kimura-book} for an introduction and detailed classification.

The infinite series of irreducible ones all have connected generic stabilizers. Using the original Sato-Kimura numbering as in \cite{kimura-book}, the basic examples are built from 
\begin{enumerate}
\item $(\SL_n, V^n)$,
\item $(\SL_n, \sym^2V^n)$,
\item $(\SL_n, \wedge^2V^n)$,
\item[(13)] $(\Sp_n, V^{2n})$,
\item[(15)] $(\OO_n, V^n)$,
\end{enumerate}
These lead to further examples by enlarging the group to contain the scalars   or replacing
$(G, V^n)$ with $(G\times \SL_m, V^n\otimes V^m)$ for certain values of $m$.

Most of the sporadic examples either have disconnected generic stabilizer
or the connectedness is not known.  A nice example is
$({\operatorname E}_6\cdot \gm, V^{27})$, which is no. 27\footnote{27=27 is a
coincidence.} on the list. The connected component of the generic stabilizer is
${\operatorname F}_4$. Since ${\operatorname F}_4$ has no outer
automorphisms, the stabilizer is ${\operatorname F}_4$, hence connected.
See also
\cite{wrightY:prehomogeneous-vector-spaces, yukie:prehomogeneous-vector-spaces, springer:some-groups-of-type-e7, katoY:rational-orbits} for several other examples.
\end{exmp}

\subsection*{Counterexamples}{\ }
\medskip

\begin{exmp}\label{reducible.exmp} Start with  $(\GL_3, \a^3_{xyz})$ and let
  $Y_1, Y_2\subset \a^3$ be cones with vertex $0$.
Then $g_1Y_1\cap g_2Y_2$ conists of
$\deg Y_1\cdot \deg Y_2$ lines  for general $g_1, g_2\in \GL_3$.
Thus, so long as $Y_1$ and $Y_2$ are not both planes, $(g_1Y_1\cap g_2Y_2)(0)=g_1Y_1\cap g_2Y_2$ is reducible and
the origin is  a good limit point. 
\end{exmp}

\begin{exmp}  Consider  $(\gm^3, \a^3_{xyz})$ and let
let $Y_1, Y_2\subset \a^3$ be surfaces that contain the $z$-axis. 
Usually $g_1Y_1\cap g_2Y_2$ is reducible, having both moving and fixed 
irreducible components for $g_1, g_2\in \gm^3$.  

If $Y_1, Y_2$ intersect transversally at the origin then the origin is not a good limit
point, even though  $g_1Y_1\cap g_2Y_2$ may have non-empty intersection with the
dense open orbit. For example, this happens for
$Y_1=(y=x^2), \ Y_2=(x=y^2)$. 

In any case, the origin is a good limit point if  $Y_1\cap Y_2$
has an irreducible component $Z$ that passes through the origin but is not contained in  a coordinate hyperplane.
\end{exmp}

\begin{exmp}  \label{2.lim.pts.exmp}
Suppose $k$ has characteristic $0$ (or at least has characteristic not equal to $2$).
Let $X\subset \p^4$ be the projective cone over a smooth quadric surface with vertex $x_0$.
Let $G\subset \aut(X)$ be the identity component. Then, $G$ acts  with 2 orbits:  $\{x_0\}$ and  $X\setminus \{x_0\}$.
Let $Y_1$ be a 2-plane contained in $X$ and containing $x_0$. We claim that
there is  a divisor $Y_2$ with 
$x_0\in Y_2\subset X$ such that $g_1Y_1^n\times_X g_2Y_2^n$ is the union of a
curve and a point for general $(g_1, g_2)\in G\times G$. Thus one preimage of $x_0$ is a good limit point, the other is not.

The computation is local at $x_0$, thus we choose affine coordinates such that
$X=(xy-uv=0)$. We can then choose  $Y_1=(x=u=0)$.  We choose $Y_2$ to be the
complete intersection  $(xy-uv=4u-(x-y)^2=0)$. 
We can eliminate $u$ to get
$$
Y_2\cong (4xy-(x-y)^2v=0)\subset \a^3.
$$
This is an irreducible hypersurface, but 
$$
4xy-(x-y)^2v=\bigl(x+y+\sqrt{1+v}(x-y)\bigr)\bigl(x+y-\sqrt{1+v}(x-y)\bigr)
$$
shows that it is non-normal along $(x=y=0)$ and  $Y_2^n$ has 2 points over the origin. 
We will next show that only one of them is a good limit point.

A typical translate of $Y_1$ by $G$ is   $(x-c^{-1}v=y-cu=0)$. Add these to the  equations
$$
x+y\pm\sqrt{1+v}(x-y)=4u-(x-y)^2=0
$$
and eliminate $u$ and $v$ to get
$$
\begin{array}{l}
x+y\pm\sqrt{1+cx}(x-y)=0 \qtq{and} \hfill(\ref{2.lim.pts.exmp}.1)\\
4c^{-1}y-(x-y)^2=0.\hfill(\ref{2.lim.pts.exmp}.2)
\end{array}
$$
Here (\ref{2.lim.pts.exmp}.2) defines a curve whose tangent line at the origin is  $y=0$. With $+$ sign, (\ref{2.lim.pts.exmp}.1) defines a curve whose tangent line at the origin is $x=0$.
So we get  an isolated intersection point at the origin. Finally,
$$
y=x\frac{\sqrt{1+cx}-1}{\sqrt{1+cx}+1}
$$
satisfies both (\ref{2.lim.pts.exmp}.1) with a $+$ sign and (\ref{2.lim.pts.exmp}.2).
\end{exmp}

\begin{exmp} \label{non-exmp.1}
	Suppose $k$ has characteristic $0$. For $X=\p^n\times \p^n$, consider 
$G = \PGL_{n+1} \times \PGL_{n+1}$,
$Y_1=\p^n\times C$, and
$Y_2=\p^n\times H$, where $C\subset \p^n$ is a smooth projective curve of positive genus and $H\subset \p^n$ a  hypersurface.  
Then 
$\dim Y_1+\dim Y_2= \dim X + n$.

For general $(g_1, g_2)\in G \times G$,
$$
g_1Y_1\cap g_1Y_2=\p^n\times (\deg C \cdot \deg H\mbox{ points}).
$$
Thus, the intersection is disconnected and its
conected components are simply connected.
So, they do not contribute to the
 fundamental group of $C$
and hence $\pi_1(g_1 Y_1 \cap g_2Y_2, z_0) \to \pi_1(Y_1, y_1)$
has infinite index for any basepoint $z_0 \in g_1 Y_1 \cap g_2Y_2$.

\end{exmp}

\begin{exmp}[Variant of
	\protect{\cite[Example, p.\
634]{hansen:a-connectedness-theorem-for-flagmanifolds-and-grassmannians}}] \label{non-exmp.2}  
Suppose $k$ has characteristic $0$.
In $X=\grass(\p^1, \p^n)$, take $G = \PGL_{n+1}$ and consider 
$$
\begin{array}{rcl}
Y_1&=&(\mbox{lines through a point})\cong \p^{n-1}, \\
  Z_L&=&(\mbox{lines in a hyperplane } L\subset \p^n)\cong \grass(\p^1, \p^{n-1}), \qtq{and}\\
Y_2&=& \cup_{c\in C} Z_{L_c},
\end{array}
$$
where  $C\subset \check\p^n$ is a smooth projective curve of positive genus
parametrizing hyperplanes $\{L_c:c\in C\}$.
Then
$$
\dim Y_1+\dim Y_2= \dim X + n-2.
$$
The key property is that 
$g_1Y_1\cap g_2 Z_L$ is either empty or is isomorphic to $\p^{n-2}$.
Thus, for general $(g_1, g_2)\in G \times G$,
$g_1Y_1\cap g_2Y_2$ is reducible and its connected components are
isomorphic to $\p^{n-2}$.  So again they do not contribute to the
fundamental group of $C$.

\end{exmp}

\begin{exmp}
	\label{example:positive-characteristic}
	Let $k$ be a field of positive characteristic, take $X = \mathbb P^2$
	with the standard action of $G = \PGL_3$.
	Let
	$Y_1 \subset \mathbb P^2$ be a line and let $Y_2 \subset \mathbb P^2$ be
	the complement of a line.
	A general translate of an intersection of $Y_1$ with $Y_2$ is isomorphic
	to $\mathbb A^1$, and the map $g_1 Y_1 \cap g_2 Y_2 \to Y_2$ is identified
	with the inclusion of a line $\mathbb A^1 \to \mathbb A^2$.
	Hence, after applying an automorphism, we may assume
	it is given by $\spec k[x,y]/(y) \to \spec k[x,y]$.
	This is not surjective on fundamental groups
	because the pullback of the connected finite \'etale Artin-Schreier
	cover $W := \spec k[x,y,t]/(t^p-t-y) \to \spec k[x,y]$
	is not connected.
\end{exmp}

\subsection*{Examples with disconnected stabilizers}{\ }
\medskip

We conclude by giving examples showing that $\pi_1(g_1 Y_1 \times_X
g_2 Y_2, z_0) \to \pi_1(Y_1, y_1) \times \pi_1(Y_2, y_2)$ may be non-surjective
when $H$ is disconnected
in the setting of Theorem~\ref{dense.orb.thm}.
For the remainder of the paper, we assume $k$ has characteristic $0$.

\begin{exmp}
	\label{example:nontrivial-fundamental-group}
	If $X$ as in Theorem~\ref{dense.orb.thm}
	has nontrivial fundamental group, we can take $Y_1 = Y_2 = X$,
	and the resulting map $\pi_1(X, x_0) \to \pi_1(X,x_0) \times \pi_1(X,x_0)$
	will fail to be surjective.
	By applying Theorem ~\ref{dense.orb.thm},
	we conclude that $H$ must be disconnected, and have at least
	$\pi_1(X, x_0)$ components.

	For a concrete example of such a variety, take $X$ to be the moduli
	space parametrizing unordered pairs of distinct points in $\p^2$.
Then $X$ is homogeneous under $\PGL_3$. The open orbit in $X\times X$ is formed by
those  $(\{p_1, p_2\}, \{p'_1, p'_2\})$  for which no three points are on a line.
The stabilizer is $\z/2\times \z/2$. The space parametrizing ordered pairs of distinct points in $\p^2$ is the universal cover of $X$. Thus
$\pi_1(X)\cong \mathbb Z/2$.
\end{exmp}

We conclude by giving a somewhat more involved example where $H$ is
disconnected,
but nevertheless $\pi_1(X,x_0) = 1$.

\begin{exmp}[Constructing $X$ and $G$]
	\label{example:disconnected-stabilizers}
	Let $X$ be the moduli space of
smooth plane conics in $\mathbb P^4$ over $\c$.
The group $\GL_5$ acts transitively on $X$ via its action on $\mathbb P^4$.

Choose coordinates $x_0, \ldots, x_4$ and a reference conic  
$C_0=(x_0^2+x_1^2+x_3^2=x_3=x_4=0)$.  The stablizer of $C_0$ is the set of matrices
\begin{align*}
	\tilde G:=
\left\{	\begin{pmatrix}
			\GO_3 & *  \\
			0 & \GL_2  		
		\end{pmatrix}
\right\}
\subset \GL_5.
	\end{align*}
Thus $\tilde G$ is connected and so is $\tilde G\cap \SL_5$ as shown by the retraction
$$
\begin{pmatrix}
			A_3 & B  \\
			0 & A_2  		
		\end{pmatrix}
\mapsto 
\begin{pmatrix}
			 d^{-1}A_3 & B  \\
			0 & dA_2  		
		\end{pmatrix}
\qtq{where} d=\det A_3\cdot \det A_2.
$$
Thus 
$X\cong \SL_5/(\tilde G\cap \SL_5)$ is simply connected.
Since we prefer faithful actions, our group $G\subset \PGL_5$ is the image of $\tilde G$. 

For $Z \subset \p^4$, let $\langle Z \rangle$ denote the linear span of $Z$ in
$\p^4$.
Observe that $G$ has a dense orbit $X^\circ \subset X$, consisting of those conics $C_1$ so that $\langle C_0\rangle\cap \langle C_1\rangle$ is
a point $p$ and neither of the $C_i$ contains $p$.

If $C_1=(x_0=x_1=x_2^2+x_3^2+x_4^2=0)$ then the  stabilizer of the ordered pair
$(C_0, C_1)$ is the set of matrices
\begin{align*}
	H:=
\left\{	\begin{pmatrix}
			\OO_2 & 0 & 0   \\
                         0  & 1 & 0      \\
			0 & 0 & \OO_2  		
		\end{pmatrix}
\right\}
\subset \PGL_5.
	\end{align*}
Thus $H$ has 4 connected components, which  can be geometrically described as
follows.  Set  $p:=\langle C_0\rangle\cap \langle C_1\rangle$. From $p$ one can
draw  two distinct tangent lines to each $C_i$. Let these tangent lines be 
$\{T^0_0, T^1_0\}$ and   $\{T^0_1, T^1_1\}$. The $H$ action permutes these lines, giving  a surjection  $H\onto \z/2\times \z/2$.

\end{exmp}

\begin{exmp}[Constructing $Y_1$ and $Y_2$]\label{example:disconnected-stabilizers-2}
	Continuing Example \ref{example:disconnected-stabilizers},
we next construct the subvarieties $Y_1$ and $Y_2 \subset X$.

Choose a  point $q \in \langle C_0\rangle\setminus C_0$ and let
$Y_1:=Y_1(q)\subset X$ be the set of those conics $C$ for which
$q \in \langle C\rangle\setminus C$.
  Then $Y_1$ is a smooth, locally closed subvariety of dimension $9$ in $X$.

To construct $Y_2$, fix 2 distinct points  $x_1, x_2\in C_0$ and 
 general $2$-planes $S_i$ such that $x_i\in S_i$. 
Let $Y_2:=Y_2(S_1, S_2)\subset X$ be the set of  conics $C$ with the following
two properties.
\begin{enumerate}
\item   We have that $c_i:=C\cap S_i$ is a single point and  $c_1\neq c_2$. 
\item Let  $\tau(C)\in
	\langle C\rangle\setminus C$ denote the intersection point of the lines
	tangent to $C$ at the $c_i$. 
Then $\tau(C)\in  \langle C_0\rangle\setminus C_0$.
\end{enumerate}
Given distinct  $c_i\in S_i$ and  $p\in  \langle C_0\rangle\setminus C_0$,
the set of such conics for $c_i = C \cap S_i$ and $\tau(C)=p$ is a principal
$\gm$-bundle
whenever $c_1, c_2$, and $p$ are not collinear.
Indeed, the  conics in $\p^2$ that pass through  $(1:0:0), (0:1:0)$
which have tangent lines at those points intersecting at $(0:0:1)$  are
precisely the hyperbolas
$x_1x_2=\lambda x_0^2$ for $\lambda \in \gm(\c)$.
Thus $Y_2$ is a smooth, locally closed subvariety of dimension $7$ in $X$.

The intersection $Z:=Y_1\cap Y_2$ consists of those conics in $Y_2$ for which
$\tau(C)=q$.  Thus  $Z$ is a smooth, locally closed subvariety of dimension $5$ in $X$.
\end{exmp}

\begin{prop}	\label{proposition:index-4}
We use the notation  of Examples~\ref{example:disconnected-stabilizers}--\ref{example:disconnected-stabilizers-2}.
	For general $(g_1, g_2) \in G \times G$, 
	$g_1 Y_1 \cap g_2 Y_2$ is irreducible and the map
	$$
		\pi_1\bigl(g_1 Y_1 \cap g_2Y_2, [C_0]\bigr) \to 
\pi_1\bigl(Y_1, [C_0]\bigr) \times \pi_1\bigl(Y_2, [C_0]\bigr)
\qtq{has index 4.}
	$$
\end{prop}

Proof.
Note first that the group action sends 
$Y_1(q)$ to $Y_1(g_1q)$ and 
$Y_2(S_1, S_2)$ to $Y_2(g_2S_1, g_2S_2)$. Thus, letting $Z = Y_1 \cap Y_2$, it is enough to show that
$$
		\pi_1\bigl(Z, [C_0]\bigr) \to 
\pi_1\bigl(Y_1, [C_0]\bigr) \times \pi_1\bigl(Y_2, [C_0]\bigr)
\qtq{has index 4.}
$$
The index is at most $4$ by Theorem \ref{dense.orb.thm} because $H$ has $4$
connected components, as shown in Example
\ref{example:disconnected-stabilizers}.
We now show the index is at least $4$.
Note that $Y_1$ has a connected degree 2 finite \'etale cover
$\widetilde{Y}_1\to Y_1$ parametrizing
pairs  $(C, c)$ where $C\in Y_1$ and $c\in C$ is one of the 2 points of $C$ whose tangent line passes through $q$.
Similarly, $Y_2$ has a connected degree $2$ finite \'etale cover
parametrizing pairs $(C,d)$ where $d$ is one of the two points of $C_0$ whose
tangent line passes through $\tau(C)$.
Let $q_1$ and $q_2$ denote the two points of $C_0$ whose tangent lines pass
through $q$.
Then, the restriction of the cover $\widetilde{Y}_1 \times
\widetilde{Y}_2$ to $Z$ splits into the 4 connected components
\[
\pushQED{\qed} 
\widetilde Z_{i,j}:=\{(C, C \cap S_i ,q_j)\} \subset Z \times \times_{Y_1 \times
Y_2} (\widetilde{Y}_1 \times
\widetilde{Y}_2)  \qtq{ for }1\leq i,j\leq 2. \qedhere
\popQED
\]

\subsection*{Conflict of Interest Statement}{\ }
On behalf of all authors, the corresponding author states that there is no conflict of interest.

\bibliography{f-h-refs}

\bigskip

  Princeton University, Princeton NJ 08544-1000, \

\email{kollar@math.princeton.edu}

\medskip

  Stanford University, Stanford CA 94305, \

\email{aaronlandesman@stanford.edu}

\end{document}